\documentclass[12pt]{amsart}

\usepackage[colorlinks]{hyperref}
\usepackage{enumerate, amssymb}

\theoremstyle{plain}
\input{diagrams.tex}
\diagramstyle[scriptlabels,height=8mm,width=8mm]

\def\bdi{\begin{diagram}}
\def\edi{\end{diagram}}

\newtheorem{thm}{Theorem}[section]
\newtheorem{cor}[thm]{Corollary}
\newtheorem{lem}[thm]{Lemma}
\newtheorem{prop}[thm]{Proposition}

\theoremstyle{definition}
\newtheorem{defi}[thm]{Definition}
\newtheorem{defis}[thm]{Definitions}
\newtheorem{conj}[thm]{Conjecture}
\newtheorem{conv}[thm]{Convention}
\newtheorem{nota}[thm]{Notation}
\newtheorem{rem}[thm]{Remark}
\newtheorem{rems}[thm]{Remarks}
\newtheorem{exa}[thm]{Example}
\newtheorem{exas}[thm]{Examples}
\newtheorem{sit}[thm]{}

\newcommand{\brem}{\begin{rem}}
\newcommand{\brems}{\begin{rems}}
\newcommand{\erem}{\end{rem}}
\newcommand{\erems}{\end{rems}}
\newcommand{\bexa}{\begin{exa}}
\newcommand{\bexas}{\begin{exas}}
\newcommand{\eexa}{\end{exa}}
\newcommand{\eexas}{\end{exas}}
\newcommand{\bdefi}{\begin{defi}}
\newcommand{\edefi}{\end{defi}}
\newcommand{\bdefis}{\begin{defis}}
\newcommand{\edefis}{\end{defis}}
\newcommand{\bcor}{\begin{cor}}
\newcommand{\ecor}{\end{cor}}
\newcommand{\blem}{\begin{lem}}
\newcommand{\elem}{\end{lem}}
\newcommand{\bconv}{\begin{conv}}
\newcommand{\econv}{\end{conv}}
\newcommand{\bconj}{\begin{conj}}
\newcommand{\econj}{\end{conj}}
\newcommand{\bprop}{\begin{prop}}
\newcommand{\eprop}{\end{prop}}
\newcommand{\bthm}{\begin{thm}}
\newcommand{\ethm}{\end{thm}}
\newcommand{\bnota}{\begin{nota}}
\newcommand{\enota}{\end{nota}}
\newcommand{\bsit}{\begin{sit}}
\newcommand{\esit}{\end{sit}}
\newcommand{\be}{\begin{eqnarray}}
\newcommand{\ee}{\end{eqnarray}}
\newcommand{\bproof}{\begin{proof}}
\newcommand{\eproof}{\end{proof}}
\def\ba{\begin{array}}
\def\ea{\end{array}}
\def\edi{\end{diagram}}
\def\bdi{\begin{diagram}}

\newcommand{\Spec}{\operatorname{Spec}}

\newcommand{\id}{\operatorname{id}}

\def\bk{{\bf k}}
\def\bK{{\bf K}}
\def\bL{{\bf L}}

\def\cJ{{\mathcal J}}
\def\cK{{\mathcal K}}

\def\cM{{\mathcal M}}

\def\cP{{\mathcal P}}
\def\cQ{{\mathcal Q}}
\def\cS{{\mathcal S}}
\def\cR{{\mathcal R}}
\def\cT{{\mathcal T}}

\def\cW{{\mathcal W}}
\def\cZ{{\mathcal Z}}

\def\tB{{\tilde B}}

\def\tB{{\tilde B}}
\def\tV{{\tilde V}}
\def\tX{{\tilde X}}

\def\tZ{{\tilde Z}}
\def\tV{{\tilde V}}

\newcommand{\A}{{\mathbb A}}

\newcommand{\C}{{\mathbb C}}
\newcommand{\Q}{{\mathbb Q}}

\newcommand{\tU}{{\widetilde U}}

\title{The Kraft-Russell Generic Equivalence Theorem and its application}

\author{Shulim Kaliman}
\address{Department of Mathematics,
University of Miami, Coral Gables, FL  33124, U.S.A.}
\email{kaliman@@math.miami.edu}

\begin{document}

\maketitle

\begin{abstract} We find some extensions of the Kraft-Russell Generic Equivalence Theorem
and using it we obtain a simple proof of a result of Dubouloz and Kishimoto.

\end{abstract}

\section{Introduction}

H. Kraft and P. Russell proved the following Generic Equivalence Theorem in \cite{KrRu}.\\

\bthm\label{in.t1}  Let $\bk$ be a field and let $p: S \to Y$ and $q: T \to Y$ be two morphisms 
of $\bk$-varieties. Suppose that 

{\rm (a)} $\bk$ is algebraically closed and of infinite transcendence degree over the prime field;

{\rm (b)} for all $y \in Y$ the two (schematic) fibers $S_y := p^{-1}(y)$ and 
$T_y := q^{-1}(y)$ are isomorphic; and

{\rm (c)} the morphisms $p$ and $q$ are affine.

Then there is a dominant morphism of finite degree 
$\varphi: X \to Y$ and an isomorphism $S \times_Y X= T \times_Y X$ over X.

\ethm

The aim of this note is to establish the following facts:

$\bullet$ the assumption (c) is unnecessary;

$\bullet$ the conclusion of Theorem \ref{in.t1} remains valid if the assumption (c) is removed
and (a) and (b) are replaced by
the following assumptions (a1) $\bk$ is an uncountable  (but not
necessarily algebraically closed) field,  (b1) there is a countable intersection $W$ of
Zariski open dense subsets of $Y$ such that $S_y $ and 
$T_y $ are isomorphic for every $y \in W$; 


$\bullet$ the conclusion of Theorem \ref{in.t1} remains valid if
 (a) and (c) are replaced by the assumptions (a2) $\bk$ is an algebraically closed
 field of finite transcendence degree over $\Q$ and 
(c2)  $p$ and $q$ are proper morphism.

Furthermore,  using Minimal Model Program over non closed  fields, Dubouloz and Kishimoto 
proved the following result \cite{DuKi}.

\bthm\label{in.t2} Let $\bk$ be an uncountable  field of characteristic zero 
and let $f : X\to S$ be a dominant morphism between
geometrically integral algebraic $\bk$-varieties. Suppose that for general closed points 
$s \in S$, the  fiber $X_s$ contains
an $\A^1$-cylinder $U_s \simeq Z_s \times \A^1$ over a 
$\kappa (s)$-variety $Z_s$. Then there exists an \'etale morphism $T\to S$ such that
$X_T=X\times_ST$ contains an $\A^1$-cylinder $U\simeq Z\times \A^1$ over a $T$-variety $Z$.

\ethm

We show by much simpler means that in the case, when $\bk$ is
an algebraically closed field (of any characteristic)  with an infinite transcendence degree over the 
prime field, the Dubouloz-Kishimoto theorem is a simple consequence of the Kraft-Russell theorem.\footnote{
Dubouloz informed the author that he and Kishimoto knew that
this version of their theorem  can be extracted from the Kraft-Russell theorem.}

\section{Assumption (c)}

The main result of this section (Theorem \ref{ac.t1}) is a straightforward adjustment of the argument  in
\cite{KrRu} (known to Russell and Kraft) but we provide it for convenience of readers.

\bnota\label{ac.n1} We suppose that $\rho : X \to Y$ is a dominant  morphism
of  algebraic $\bk$-varieties where $\bk$ is an algebraically closed field with
an infinite transcendence degree over its prime field. 
Recall that 
there is a field $\bk_0\subset \bk$ which is finitely generated over the prime field and 
a morphism $\rho_0 : X_0\to  Y_0$ of $\bk_0$-varieties such that the morphism $\rho : X \to Y$
is obtained from $\rho_0$ by the base extension $\Spec \bk \to \Spec \bk_0$.
Denote by $\bK_0$ the field of rational functions on $Y_0$, i.e. $\Spec \bK_0$
is the generic point of $Y_0$. Put $X_{0,\omega} = X_0 \times_{\Spec \bk_0} \Spec \bK_0$.
\enota

The next fact was proven in \cite[Lemma 1]{KrRu}) (but unfortunately under 
the additional unnecessary assumption that $X$ is affine).

\blem\label{ac.l1} Let Notation \ref{ac.n1} hold. Then every
$\bk_0$-embedding $\bK_0 \hookrightarrow \bk$ defines a closed point $y \in Y$
and an isomorphism 
$$X_{0,\omega} \times_{\Spec \bK_0} \Spec \bk \to X_y=\rho^{-1} (y).$$

\elem 

\bproof  Without loss of generality we suppose that $Y$ is affine.
Let $\bk_0 [Y_0]$ be the algebra of regular functions on $Y$. Following \cite[Lemma 1]{KrRu} we see that since $\bk_0 [Y_0] \subset \bK_0$ any $\bk_0$-embedding
$\bK_0 \hookrightarrow \bk$ yields a $\bk_0$-homomorphism
$\bk [Y]=\bk_0 [Y_0]\otimes_{\bk_0} \bk \to \bk$ and, thus, a closed point $y$ in $Y$. 
Let $U_0$ be a Zariski dense open affine subset of $X_0$, $U_{0,\omega} = U_0 \times_{\Spec \bk_0} \Spec \bK_0$,
$U= U_0 \times_{\Spec \bk_0} \Spec \bk$ and $U_y$ be the fiber over $y$ of the restriction $U \to Y$ of $\rho$. 
Continuing the argument of Kraft and Russell we have
\be\label{ac.eq1}   U_{0,\omega}  \times_{\Spec \bK_0} \Spec \bk \simeq U_0  \times_{Y_0} \Spec \bk \simeq U  \times_{Y} \Spec \bk  = U_y.                            \ee
Furthermore, if $V_0$ is a Zariski open affine subset of $U_0$  then the way the isomorphism 
$U_{0,\omega}  \times_{\Spec \bK_0} \Spec \bk \simeq  U_y$ was constructed in Formula (\ref{ac.eq1}) yields the
commutative diagram
\be\label{ac.eq0} \begin{array}{ccc} V_{0,\omega}  \times_{\Spec \bK_0} \Spec \bk & \,\,\simeq  & \, \,  V_y \\   
\Big\downarrow & & \Big \downarrow\\
U_{0,\omega}  \times_{\Spec \bK_0} \Spec \bk &\,\, \simeq  & \,\,U_y  \end{array} \ee 
where the vertical arrows are the natural embeddings (in other works, one has an isomorphism between
the structure sheaves of  $U_{0,\omega}  \times_{\Spec \bK_0} \Spec \bk$ and $U_y$).  
Consider a covering of $X_0$  (resp. $X_{0,\omega}$, resp. $X_y$) by affine charts $\{ U_0^i \}_{i=1}^n$
(resp. $\{ U_{0,\omega}^i \}_{i=1}^n$, resp. $\{ U_y^i \}_{i=1}^n$). If $n=2$ then 
applying Diagram (\ref{ac.eq0}) for the embeddings
$U_{0,\omega}^1\cap U_{0,\omega}^2\hookrightarrow U_{0,\omega}^i$ and $U_y^1\cap U_y^2\hookrightarrow U_y^i$ and gluing the affine charts
 we get an isomorphism between
$X_{0,\omega} \times_{\Spec \bK_0} \Spec \bk$ and $X_y=\rho^{-1} (y)$.  Furthermore, we see that Diagram (\ref{ac.eq0}) remains true when
$U$ (resp. $V$) is not affine but only a union of two affine sets. Then the similar argument and the induction by $n$ yields
the desired isomorphism
$X_{0,\omega} \times_{\Spec \bK_0} \Spec \bk \to X_y$ for $n \geq 3$.
\eproof

\bnota\label{ac.n2}
Let $\varphi : Z \to X$ be a morphism of
 algebraic $\bk$-varieties and $\bk_0$ be a subfield of  $\bk$ such that for some
$\bk_0$-varieties $X_0$ and $Z_0$ one has 
$X=X_0\times_{\Spec \bk_0} \Spec \bk$ and $Z=Z_0\times_{\Spec \bk_0} \Spec \bk$.
 Suppose that $\bk_1 \subset \bk$ is a finitely generated extension 
of $\bk_0$ such that for $\bk_1$-varieties $X_1$ and $Z_1$
 there exists a morphism $\varphi_1 : Z_1\to X_1$  for which $\varphi$
is obtained from $\varphi_1$ by the base extension $\Spec \bk \to \Spec \bk_1$.
However, besides $\bk_0$ the description of $\varphi$ requires not the whole field $\bk_1$
but only a finite number of elements of $\bk_1$ (because $\varphi$ is defined by the
homomorphisms of rings of regular functions on affine charts and these rings are
finitely generated). Thus for the $\bk_0$-algebra $C\subset \bk_1$ generated by 
these elements we have the following observation used by Kraft and Russell in their proof
for the affine case.
\enota

\blem\label{ac.l3} Let $X$ and $Z$ be algebraic varieties over a field $\bk$
and $\varphi : Z \to X$ be a morphism. Suppose that $\bk_0$, $X_0$ and $Z_0$
are as before.
Then there exist  a finitely generated $\bk_0$-algebra
$C\subset \bk$, ringed spaces $\tX$ and $\tZ$ with structure
sheaves consisting of $C$-rings\footnote{
If $C$ is a subring of a ring $R$ we call $R$ a $C$-ring and a homomorphism
of two $C$-rings whose restriction to $C$ is the identity map is called a $C$-homomorphism.} 
and a $C$-morphism $\tilde \varphi : \tZ \to \tX$
such that $X=\tX\times_{\Spec C} \Spec \bk$, $Z=\tZ\times_{\Spec C} \Spec \bk$
and $\varphi = \tilde \varphi \times_{\Spec C} \id_{\Spec \bk}$.

\elem

\bthm\label{ac.t1} The Generic Equivalence Theorem is valid without the assumption (c).

\ethm

\bproof  As before we can choose a field $\bk_0\subset \bk$ which is finitely generated over the prime field 
such that for some morphisms $p_0: S_0 \to Y_0$ and $q_0: T_0 \to Y_0$ of $\bk_0$-varieties
the morphisms $p : S \to Y$ and $q: T \to Y$ are obtained from these ones
via the base extension $\Spec \bk \to \Spec \bk_0$. Suppose that $\bK_0$ is the field
of rational functions on $Y_0$. 

As in \cite{KrRu} by Lemma \ref{ac.l1} we get the following isomorphisms in self-evident notations
$$S_{0,\omega} \times_{\Spec \bK_0} \Spec \bk \simeq S_y \simeq T_y \simeq 
T_{0,\omega} \times_{\Spec \bK_0} \Spec \bk.$$ 
By Lemma \ref{ac.l3}, for the isomorphism $S_{0,\omega} \times_{\Spec \bK_0} \Spec \bk \simeq 
T_{0,\omega} \times_{\Spec \bK_0} \Spec \bk$ there exists a finitely generated
$\bK_0$-algebra $C$ in $\bk$ such that one has 
$$S_{0,\omega} \times_{\Spec \bK_0} \Spec C \simeq 
T_{0,\omega} \times_{\Spec \bK_0} \Spec C.$$
 Choosing a maximal ideal $\mu$ of $C$
contained in the image of the morphism $S_{0,\omega} \times_{\Spec K_0} \Spec  C \to \Spec C$
and letting $\bL_0=C/\mu$ we get 
\be\label{ac.eq2} 
S_{0,\omega} \times_{\Spec \bK_0} \Spec \bL_0 \simeq 
T_{0,\omega} \times_{\Spec \bK_0} \Spec \bL_0.\ee
By construction the field $\bL_0$ is a finite
extension of $\bK_0$. It follows that there is a finite extension $\bL$ of the field $\bK$
of rational functions on $Y$ such that
$$S_\omega \times_{\Spec \bK} \Spec \bL \simeq T_\omega \times_{\Spec \bK} \Spec \bL$$
where $S_\omega$ and $T_\omega$ are generic fibers of $p$ and $q$ respectively.
Since  $S_\omega \times_{\Spec \bK} \Spec \bL \simeq S \times_{Y} \Spec \bL$
and $T_\omega \times_{\Spec \bK} \Spec \bL \simeq T \times_{Y} \Spec \bL$ there is a
dominant morphism $X \to Y$ for which $S\times_Y X \simeq T\times_YX$ and we are done.
\eproof

\brem\label{ac.r2} It is interesting to discuss what happens to Theorem \ref{ac.t1} 
if the field $\bk$ is not algebraically closed (but still of infinite transcendence degree over the prime field). Then there may be no embedding $\bK_0 \hookrightarrow \bk$
as in Lemma \ref{ac.l1}.  However, for a finite extension $\bk_1$ of $\bk$ one can find an embedding $\bK_0 \hookrightarrow \bk_1$.
Consider the morphisms 
$p_1 : S_1 \to Y_1$ and $q_1 : T_1 \to Y_1$ of $\bk_1$-varieties obtained from
$p: S \to Y$ and $q: T \to Y$ via the base extension $\Spec \bk_1 \to \Spec \bk$. 
Then until Formula (\ref{ac.eq2}) the argument remains valid with $\bk$, $p$ and $q$ replaced by $\bk_1$, $p_1$ and $q_1$.
In Formula (\ref{ac.eq2}) the field $\bL_0$
may contain a nontrivial finite extension $\bk_0^1$ of $\bk_0$. Taking a bigger field $\bk_1$ we can suppose that $\bk_0^1$ is a subfield of $\bk_1$
and proceed with the proof.
Hence, though we cannot get the exact formulation of the Generic equivalence theorem in the case of non-closed fields,
we can claim that for a finite extension $\bk_1$ of $\bk$ and $S_1,T_1$ and $Y_1$ as before there is a dominant morphism 
of $\bk_1$-varieties of finite degree 
$X_1\to Y_1$ and an isomorphism $S_1 \times_{Y_1} X_1\simeq T_1 \times_{Y_1} X_1$ over $X_1$.
\erem

\section{Very general fibers and non-closed fields}

It is obvious that the assumption that an isomorphism $S_y \simeq T_y$ holds for every 
$y \in Y$ in the Kraft-Russell theorem can be replaced with the assumption that it
is true for a general point of $Y$, i.e. for every point contained in some Zariski open dense
subset of $Y$. However, the author does not know whether the proof of Kraft and Russell
can be adjusted to the case when $y$ is only a very general point of $Y$, i.e.
it is in a complement of the countable union of proper closed subvarieties of $Y$.
Hence we shall use a different approach. Namely, we shall use the technique
which was communicated to the author by Vladimir Lin in 1980s and which
was used in his unpublished work with Zaidenberg on a special case of The
Generic Equivalence Theorem. The negative feature of this new proof is that
we have to work over an {\bf uncountable} field $\bk$. 
However, we do not require that this field is
algebraically closed.

\bdefi\label{abc.d1} We say that an uncountable subset $W$ of an algebraic $\bk$-variety
$X$ is Zariski locally dense if $W$ is not contained in any countable union of proper closed suvarieties of $X$.
\edefi

\bexa\label{abc.exa1} Let $W$ be the complement of a countable union $\bigcup_{i=1}^\infty Y_i$
of closed proper subvarieties of $X$. Then $W$ is Zariski locally dense.
Indeed, assume the contrary. That is, $W$ is contained in a union $\bigcup_{i=1}^\infty Z_i$
of proper closed subvarieties of $X$ and $X=\bigcup_{i=1}^\infty Y_i \cup \bigcup_{i=1}^\infty Z_i$. Without loss of generality we can suppose that $X$
is affine and using a finite morphism of $X$ onto some affine space $\A_\bk^n$ we reduce
the consideration to the case of $X \simeq \A_\bk^n$. Note that equations  of all $Y_i$'s and $Z_i$'s
involve a countable number of coefficients. Let $\bk_0$ be the smallest subfield of $\bk$
containing all these coefficients.
Since $\bk_0$ is countable there are points in $\A_\bk^n$ whose coordinates are algebraically
independent over $\bk_0$. Such a point cannot be
contained in $\bigcup_{i=1}^\infty Y_i \cup \bigcup_{i=1}^\infty Z_i$. A contradiction.
\eexa

The aim of this section is the following.

\bthm\label{abc.t1} Let  $\bk$ be an uncountable field of characteristic zero and $p : S \to Y$ and $q: T\to Y$
be morphisms
of  $\bk$-varieties. Suppose that $W$ is a Zariski locally 
dense subset of $Y$ and for every $y \in W$ there is an isomorphism
$p^{-1}(y)=S_y \simeq T_y=q^{-1}(y)$. Then there is a dominant morphism of finite degree
$X\to Y$ such that $S\times_Y X$ and $T\times_Y X$ are isomorphic over $X$.

\ethm

The proof requires some preparations.
We start with the following simple fact.

\bprop\label{abc.p1} Let $Y$ be an algebraic $\bk$-variety, 
 $X$ and $Z$ be  subvarieties of $Y\times \A_\bk^n$,  $\rho : X \to Y$ and 
 $\tau: Z \to Y$ be  the natural projections, and 
 $P$ be an algebraic family of rational maps $\A_\bk^n \dashrightarrow \A_\bk^n$.
 Suppose that $\cP$ is a subvariety of $Y\times P$ such that for every $(y,f) \in \cP$
 the map $f$ is regular on $X_y = \rho^{-1}(y)$.
Let  $\cP_{X,Z}$ be the subset of $\cP$ that 
consists of all elements $(y,f)$ such that $f(X_y) \subset Z_y$ for $Z_y:=\tau^{-1}(y)$.
Then $\cP_{X,Z}$ is a constructible set.

\eprop

\bproof 
Consider the morphism $\kappa: X \times_Y \cP\to Y\times \A_\bk^1$ given
by $(x, f) \to (\rho (x),f(x))$. Then $(Y\times \A^n_\bk) \setminus Z$ and, therefore,
$\kappa^{-1}((Y\times \A^n_\bk) \setminus Z))$
are constructible sets.
The image $R$ of the latter under the natural projection $ X \times_Y \cP \to  \cP$
is a constructible set by the Chevalley's theorem (EGA IV, 1.8.4). Note that $\cP\setminus R$ 
coincides with $\cP_{X,Z}(N)$, i.e. it is also constructible and we are done. 
\eproof

Letting $Z=Y \times o$ where $o$ is the origin in $\A_\bk^n$ we get the following.

\bcor\label{abc.c1} The subset $\cP_X^0(N)$ of $\cP$ that consists of all elements $(y,f)$ such that
$f$ vanishes on $X_y$ is a constructible set.

\ecor

\bnota\label{abc.n1} 
(1) Let  $P (N)$ 
 consist of $2n$-tuples $\varphi =(f_1,g_1,f_2,g_2, \ldots, f_n,g_n)$ of polynomials
on $\A_\bk^n$ of degree at most $N$ such that $g_1, \ldots , g_n$ are not zero polynomials.
We assign to $\varphi$
the rational map $ \breve \varphi : \A_\bk^n \dashrightarrow \A_\bk^n$ given by
$ \breve \varphi=(\frac{f_1}{g_1}, \ldots, \frac{f_n}{g_n}) $ and denote the variety 
of such rational maps by $R(N)$ with $\Theta : P(N) \to R(N)$  being the morphism  given by 
$\Theta (\varphi )=\breve \varphi$.

(2) Let the assumptions of Theorem \ref{abc.t1} hold and $Y$ be affine. 
Consider a cover of $S$ (resp. $T$) by  a collection $\cS=\{S^i\}_{i\in J}$ of  affine charts
(resp. $\cT=\{T^i\}_{i\in J}$) where $J$ is finite set of indices.
We can suppose  that for some $n>0$
every $S^i$ (resp. $T^i$) is a closed subvariety of $Y\times A_{\bk}^n$
where the natural projection $S^i \to Y$ is the restriction of $p$ (resp. $T_i \to Y$ is the restriction of $q$).
By  $S_y^i$ (resp. $ T_y^i$) we denote  $S_y \cap S^i$ (resp. $T^i\cap T_y$).
We treat each transition isomorphism $\alpha^{ij}$ of $\cS$ (resp. $\beta^{ij}$ of $\cT$)
as the restriction of some rational map $\A_Y^n \dashrightarrow \A_Y^n$ and, choosing $N$ large
enough we assume that for every $y \in Y$ each  rational map  $\alpha^{ij}|_{S_y^i}$ (resp. $\beta^{ij}|_{T_y^i}$)
is contained in $\Theta (P(N))$.

(3)
Suppose that $\cQ(N)=\prod_{i,j\in J} P(N)$,
i.e. each element of $\cQ (N)$ is $\Phi= \{ \varphi^{ij}\in P(N)| i,j \in J \}$
where $\varphi^{ij} =(f_1^{ij},g_1^{ij},f_2^{ij},g_2^{ij}, \ldots, f_n^{ij},g_n^{ij})$.
Let $F(N)$ be the subset of $Y\times \cQ (N)$ consisting of all elements $(y, \Phi)$
such that for all  $i,j,i',j'\in J$ and $y \in Y$ one has the following
\be\label{abc.eq2} \breve \varphi^{i'j'} \circ \alpha^{ii'}|_{S_y^i}  =  
\beta^{jj'} \circ \breve \varphi^{ij}|_{S_y^i};\ee  \be\label{abc.eq3}    \forall x \in S^i \,\, \exists j \text{ such that } \forall k=1, \ldots , n \,\, g_k^{ij} (x) \ne    0;                     \ee
\be\label{abc.eq4} 
\breve \varphi^{ij} (S_y^i)\subset T_y^j. \ee
\enota

\blem\label{abc.l1} The set $F (N)$ is constructible.

\elem

\bproof   

Every coordinate function of the rational map
$$(\breve \varphi^{i'j'} \circ \alpha^{ii'} -  \beta^{jj'} \circ \breve \varphi^{ij}): \A_\bk^n\dashrightarrow \A_\bk^n $$
can be presented as a quotient of two polynomials with the degrees of the numerator
and the denominator  bounded by  some constant $M$ depending on $N$ only.
That is, the ordered collection $\nu_{i,j,i',j'}$  of the numerators 
of this rational map is contained in $P(M)$.
Consider the morphism
$\tilde \nu_{i,j,i',j'} :  Y\times \cQ (N) \to Y  \times P(M)$
where $\tilde \nu_{i,j,i',j'} = ({\rm id}, \nu_{i,j,i',j'})$.
Let $\cZ_i$  be the subvariety of $Y  \times P(M)$ that consists of all elements
$(y,f_1, \ldots , f_n)$ such that $f_k|_{S_y^i}\equiv 0$ for every $k$. 
By Corollary \ref{abc.c1} $\cZ_i$ is a constructible set.
Hence its preimages $\tilde \cZ_{i,j,i',j'}$ in $Y\times \cQ (N)$
under $\tilde \nu_{i,i',j,j'} $ is also constructible.
Note that  the variety
$C=\bigcap_{i,j,i',j' \in J} \tilde \cZ_{i,j,i',j'}$ consists of all elements satisfying
Formula (\ref{abc.eq2}).

Consider the morphism $\theta_{ij} : S^i\times_Y C \to S^i \times \A_\bk^n$ over $S^i$
which sends each point $(x,\Phi)$ to $(x, g_1^{ij} (x), \ldots ,  g_n^{ij} (x))$.
Let $L$ be the union of the coordinate hyperplanes in $\A_\bk^n$ and $L_{ij}=\theta_{ij}^{-1} (S^i \times L)$ .
Then $K_i=\bigcap_{j\in J}L_{ij}$ is the subvariety of $S^i\times_Y C$ 
that consists of all points $(x,\Phi)$ such that for every  $j \in J$  there exists $1\leq k \leq n$
with $g_k^{ij} (x)= 0$. Let $\cK_{i}$ be the image of $K_{i}$ in $C$
under the natural projection $S_{i} \times_Y C \to C$ i.e. it  is constructible by the Chevalley's theorem. 
Then its complement 
$\cM_i$ consists of elements $(y, \Phi)\in C$ such that for every $x \in S_y^i$ there exists $j \in J$
for which $g_k^{ij} (x)\ne 0$ for every $k=1, \ldots, n$.
Hence the constructible set $D =\bigcap_{i\in J}  \cM_{i}$ satisfies Formula (\ref{abc.eq3}).

Let $R^{ij}$ be the subvariety of $S^i  \times_Y D$ that consists of points $(x,\Phi)$ for which
$g_k^{ij} (x)\ne 0$ for every $k=1, \ldots, n$. That is, the map $\kappa_{ij}:
R^{ij} \to Y \times \A_\bk^n$ sending each point $(x,\Phi)$ to $(p(x), \breve \varphi^{ij} (x))$ is
regular.  Let $\cR^{ij}\subset R^{ij}$ be the preimage of $\A_\bk^n\setminus T^j$ 
under this map. Then $\cR_i = \bigcup_{j \in \cJ} \cR^{ij}$ is a constructible subset of $S^i  \times_Y D$
and, therefore, its image $R_i$ in $D$ under the natural projection $S^i\times_YD \to D$
is also constructible. Note that $F(N):=D \setminus \bigcup_{i \in J} R_i$ satisfies Formula (\ref{abc.eq4})
and we are done. \eproof

\brem\label{abc.r1}  Formulas (\ref{abc.eq2}), (\ref{abc.eq3}) and (\ref{abc.eq4})
guarantee that any point $(y,\{ \varphi^{ij} | i,j \in J\})$ in $F(N)$ defines a morphism
$S_y \to T_y$. Hence we treat  $F(N)$ further as
collections of such morphisms.

\erem

\bnota\label{abc.n2} Exchanging the role of $S$ and $T$ we get a constructible set $G(N)$, i.e.
each element of $G$ defines a morphism
$T_y \to S_y$. In particular,  $F(N)\times_YG(N)$ consists of elements
$\{ (y,f_y,g_y)\}$ where $f_y:S_y\to T_y$ and $g_y: T_y \to S_y$ are morphisms.

\enota

\blem\label{abc.l2} Let $H(N)$ be the subset of $F(N)\times_YG(N)$ that consists of all elements
$(y,f_y,g_y)$
such that each $f_y$ is an isomorphism and $g_y=f_y^{-1}$. Then $H(N)$ is a constructible
set.

\elem

\bproof Note that each element $h=(y,f_y,g_y)$ of $F(N)\times_Y G(N)$ defines the morphism
$\kappa_h : S \to S\times_YS$ which sends $s \in S_y$ to $(s, g_y \circ f_y (s))$. Let $\Delta_S$ be
the diagonal in $S\times_YS$ and let $H'(N) \subset F(N)\times_Y G(N)$ be the subset
that consists of those elements $h$ for which $\kappa_h(S) \subset \Delta_S$. By Proposition \ref{abc.p1}
$H' (N)$ is constructible. Exchanging the role of $S$ and $T$ we get the similar
constructible set $H''(N)$. Letting $H(N)=H'(N) \cap H''(N)$ we get the desired conclusion.

\eproof

\blem\label{abc.l3}  Let the assumptions of Theorem \ref{abc.t1} hold, 
$H(N)$ be as in Lemma \ref{abc.l2} and
$W(N)=\rho (H(N))$ where $\rho : H(N) \to Y$ is the natural projection.
Then for some number $N$ the set $W(N)$   contains a Zariski dense open subset of $Y$.
\elem

\bproof   Note  that for any isomorphism $\varphi : S_y \to T_y$ 
we can find $N$ for which $(\varphi, \varphi^{-1})$ is an element of $H(N)$.
Hence the  assumptions of Theorem \ref{abc.t1} imply that $Y=\bigcup_{N=1}^\infty W(N)$.
Therefore, one of $W(N)$'s is Zariski locally dense in $Y$.
Furthermore, it is constructible by the  Chevalley's theorem which implies that it
contains a Zariski open dense subset $U$ of $Y$.  This is the desired conclusion. 
\eproof

\noindent {\em Proof of Theorem \ref{abc.t1}}. Without loss of generality we
suppose that $Y$ is affine.
Let $N$ be as in Lemma \ref{abc.l3}, $H=H(N)$ and $\rho : H\to Y$ be the natural morphism.
It is dominant by Lemma \ref{abc.l3}. Taking a smaller $H$ we can suppose that
it is affine. Then we have the natural embedding $\rho^* : \bk [Y] \hookrightarrow \bk [H]$
of the rings of regular functions. For the field $K$ of rational functions on $Y$
consider the $K$-algebra $A=K\otimes_{\bk [Y]} \bk [H]$. By the Noether normalization lemma 
one can find 
algebraically independent elements $z_1, \ldots , z_k \in \bk [H]$ such that $A$
is a finitely generated over the polynomial ring $K[z_1, \ldots , z_n]$. 
Choose elements
$b_1, \ldots , b_k\in \bk$ so that the subvariety $X$ of $H$ given by
the system of equations $z_1-b_1= \ldots = z_k-b_k=0$ is not empty.
Then the field of rational functions on $X$ is a finite extension of $K$, i.e.
we get a dominant morphism $X\to Y$ of finite degree.

Note that we can veiw $x \in X \subset H$ as an element $(y,f_y,g_y)$ of
$F(N)\times_Y G(M)$ as in Lemma \ref{abc.l2} 
such that $y=\rho (x)$ and $f_y : S_y \to T_y$ is an isomorphism while $g_y : T_y\to S_y$ is its inverse.
Hence the map $S\times_Y X \to T\times_Y X$ that sends
every point $(s, x) \in S\times_Y X$ to $(f_y(s), x)$ is an isomorphism.
This is the desired conclusion.   \hfill $\square$

\brem\label{abc.r2}
 We do not know if the morphism $X \to Y$ in Theorem \ref{abc.t1}
 can be made \'etale in the case of a positive characteristic.
However, if $\bk$ has characteristic zero then over a Zariski dense open subset $U$ of $Y$
this morphism is smooth by the Generic Smoothness theorem \cite[Chapter III, Corollary 10.7]{Har} and replacing
$Y$ by $U$ we can suppose that $X \to Y$ is \'etale.

\erem

\section{Case of $\bar \Q$-varieties}

\bnota\label{cq.n1} In this section $\bk_0$ is an algebraically closed field of finite transcendence degree
over $\Q$  (e.g., $\bk_0$ is the field $\bar \Q$ of algebraic numbers) and
$p_0 : S_0 \to Y_0$ and $q_0 : T_0 \to Y_0$ are morphisms
of algebraic $\bk_0$-varieties. By the Lefschetz principle we treat $\bk_0$ as a subfield
of $\C$ and we denote by
 $p : S \to Y$ and $q : T\to Y$ are complexifications of these morphisms $p_0$ and $q_0$
(i.e., say, $S$ coincides with $S_0\times_{\Spec \bk_0} \Spec \C$). 
\enota

The analogue of the Kraft-Russell theorem for $k_0$-varieties
can be reduced to the complex case if the following is true.

\bconj\label{cq.con1} {\em Let Notation \ref{cq.n1} hold and the fibers
$p_0^{-1}(y)$ and $q_0^{-1}(y_0)$ be isomorphic for general points $y_0\in Y_0$.
Then the fibers $S_y=p^{-1}(y)$ and $T_y=q^{-1}(y)$  are isomorphic for  general points $y \in Y$.}
\econj

We can prove this conjecture only in the case of proper morphisms
$p_0 : S_0 \to Y_0$ and $q_0 : T_0 \to Y_0$, and our proof is based on the theory
of deformations of compact complex spaces.

\bdefi\label{cq.d1}  
 A deformation of a compact complex space $Z$ is a proper flat
holomorphic map $\rho :  \cZ \to B$ of a complex spaces such that for a marked 
point $b_0\in B$ one has $ \rho^{-1}(b_0)=Z$. A deformation $ \rho$
is called versal if for any other deformation 
$\kappa : \cW \to D$ of $Z$ with $Z=\kappa^{-1}(d_0)$ there is 
 a holomorphic map $g: (D,d_0)\to (B,b_0)$ of the germs
such that $g^*( \rho)=\kappa|_{(D,d_0)}$.
\edefi

We need the following crucial results of Palamodov \cite[Theorem 5.4]{Pa76} and  \cite{Pa73}. 

\bthm\label{cq.t1} {\rm (1)} Every compact complex space $Z$ is a fiber $\rho^{-1} (b_0)$ of a proper
flat map $\rho : \cZ \to B$ of complex spaces which is a versal deformation of each of its fibers.

{\rm (2)} The space $\cM$ of classes of isomorphic complex spaces admits
a $T_0$-topology such that every proper flat family $\rho : \cZ \to B$ induces a continuous 
map $\theta : B \to \cM$. 

{\rm (3)} Furthermore, if $\rho$ in (2) is a versal deformation at every point of $B$
then $\theta$ is an open map.
\ethm

\bthm\label{cq.t2} Let Notation \ref{cq.n1} hold and  the morphisms $p$ and $q$ be proper.
Then 
Conjecture \ref{cq.con1}  is true.
\ethm

\bproof
Without loss of generality we suppose that $Y_0$ is affine, i.e. we view $Y_0$ as a
closed subvariety of $ \A_{\bk_0}^n$. Hence $Y$ is a closed subvariety of $\C^n$ and we
can treat the set of points in $Y$ whose coordinates are in $\bk_0$ as $Y_0$.
Let $Y_1$ be the closure of $Y_0$ in $Y$ in the standard topology (i.e. $Y_1$ contains all
points of $Y \subset \C^n$ with real coordinates).  
That is, $Y_1$ is Zariski locally dense
in $Y$ in the terminology of Definition \ref{abc.d1}. Hence by Theorem \ref{abc.t1} it
suffices to establish isomorphisms $S_y \simeq T_y$ for  a general $y_0 \in Y_1$.

Let $\rho : \cZ \to B$ be a versal deformation as in Theorem \ref{cq.t1} (1) for $Z=S_{y_0} \simeq T_{y_0}$
and let $\theta : B \to \cM$ be as in Theorem \ref{cq.t1} (2). 
For some 
neighborhood $Y'$ (in the standard topology) of $y_0$ in $Y$ 
we have holomorphic maps $\hat p : (Y',y_0)\to (B, b_0)$
and $\hat q : (Y',y_0)\to (B, b_0)$ such that $(\hat p)^*(\rho)=p|_{Y'}$ and $(\hat q)^*(\rho)=q|_{Y'}$.
To prove that $S_y$ and $T_y$ are isomorphic it suffices to prove that they
are biholomorphic (by virtue of [SGA 1, XII, Theorem 4.4]). That is, 
it suffices for us to establish the equality $p':=\theta \circ \hat p = \theta \circ \hat p=:q'$
and, furthermore, as we mentioned before  it is enough to establish 
equality  $p'|_{Y_1'} = q'|_{Y_1'}$ where $Y_1'=Y_1\cap Y'$.

Assume the contrary. Then by Theorem \ref{abc.t1}
the set $R_0\subset Y_1'$ of points $y$ for which $p'(y)=q'(y)$ cannot be Zariski locally dense,
i.e. it is contained in a countable union of subsets of $Y_1'$ which are nowhere dense in $Y_1'$.
Let $R_1\subset Y_1'$ (resp. $R_2\subset Y_1'$) be the set of points $y$ such that there is a neighborhood 
$U_y \subset \cM$ of $p' (y)$ that does not contain 
$q' (y)$ (resp. a neighborhood 
$V_y \subset \cM$ of $q' (y)$ that does not contain
$p' (y)$). Since $\cM$ is a $T_0$-space we see that
$R_0 \cup R_1\cup R_2=Y_1'$. Furthermore, since the map $\theta$ is open
(by Theorem \ref{cq.t1}(3)) we can suppose that $U_y=\theta (\tU_y)$ (resp. 
$V_y=\theta (\tV_y)$) where $\tU_y$ is a neighborhood of $\hat p (y)$ in $B$
(resp.  $\tV_y$ is a neighborhood of $\hat q (y)$ in $B$).  Since $B$ is a germ
of a complex space we can consider a metric on it which induces the standard topology. 
Let $R_1^n$ be the set of points
$y \in R_1$ such that $\tU_y$ contains the ball  $D(y,\frac{1}{n})\subset B$ of radius $\frac{1}{n}$ (in this metric)
with center at $y$ and let $R_2^n$ be the similar subset of $R_2$.
Then we have $$Y_1'=R_0 \cup \bigcup_{i=1}^\infty R_1^n \cup \bigcup_{i=1}^\infty R_2^n.$$
By the Baire category theorem there is a nonempty open subset $W \subset Y_1'$ 
and $n$ such that, say, $R_1^n$ is everywhere dense in $W$. In particular,
for every point $y _1\in W\cap R_1^n$ the ball $D(\hat p (y_1), \frac{1}{n})$
does not meet $\theta^{-1} ( q' (y_1))$.
Without loss of generality we can suppose that $\hat p (y_1)$ is a smooth point of
$\hat p (Y')$ and taking a larger $n$ we can also assume that 
$\hat p (W)$ coincides with the intersection of $\hat p(Y_1')$ with $D(\hat p (y_1), \frac{1}{2n})$.
Hence we can choose a point $y_2 \in W \cap Y_0$ near $y_1$
such that for $b_2:=\hat p (y_2)$ the ball  $D(b_2, \frac{1}{3n})$ does not meet $\theta^{-1} (q' (R_1^n\cap W))$.

On the other hand by the assumption of Conjecture \ref{cq.con1}
 we have $S_{y_2}\simeq T_{y_2}$ and since $\rho : \cZ \to B$
is a versal deformation at every point of $B$ (by Theorem \ref{cq.t1} (1)) there exists
a  map $\breve q : (Y',y_2)\to (B, b_2)$ such that $(\breve q)^*(\rho)=q|_{(Y', y_2)}$.
By continuity $D(b_2, \frac{1}{3n})$ contains points from $\breve q (R_1^n\cap W))$.
Thus it must contain points from $\theta^{-1} (q' (R_1^n\cap W))$
since $q'=\theta \circ \hat q$ and $\theta \circ \hat q (R_1^n\cap W)= \theta \circ \breve q (R_1^n\cap W)$.
This contradiction shows that $R_0$ is Zariski locally dense in $Y$.
Now the desired conclusion follows from Theorem \ref{abc.t1}.

\eproof

\brem\label{cq.r1} The assumption that $\bk_0$ is algebraically closed can be dropped
from the formulation of Theorem \ref{cq.t2} since it is not used in the proof.

\erem

\bthm\label{cq.t3} Let Notation \ref{cq.n1} hold and let
$p_0: S_0 \to Y_0$ and $q_0: T_0 \to Y_0$ be proper morphisms. 
Suppose that 
for all $y_0 \in Y_0$ the two (schematic) fibers $p_0^{-1}(y_0)$ and 
$q_0^{-1}(y_0)$ are isomorphic.
Then there is a dominant morphism of finite degree 
$X_0\to Y_0$ and an isomorphism $S_0 \times_{Y_0} X_0= T_0 \times_{Y_0} X_0$ over $X_0$.

\ethm

\bproof  
Let $\bK_0$ be the field of rational functions on $Y_0$ and $\omega$ be the generic point of $Y_0$.
Since $\C$ is algebraically closed and of infinite transcendence degree over $\bk_0$
we can always find an injective homomorphism $\bK_0 \hookrightarrow \C$.
Then by Lemma \ref{ac.l1} this homomorphism  defines a closed point $y \in Y$ and
$$S_{0,\omega} \times_{\Spec \bK_0} \Spec \C \to S_y$$
where $S_{0,\omega}$ is the generic fiber of $p_0$.  Since the similar fact holds
for $q_0: T_0\to Y_0$ and since $S_y \simeq T_y$ by Theorem \ref{cq.t2}  we have
$$S_{0,\omega} \times_{\Spec \bK_0} \Spec \C \simeq S_y \simeq T_y \simeq 
T_{0,\omega} \times_{\Spec \bK_0} \Spec \C.$$ 
Then repeating the argument from the proof of Theorem \ref{ac.t1} 
we construct a finite extension $\bL_0$ of $\bK_0$ for which
$$S_{0,\omega} \times_{\Spec \bK_0} \Spec \bL_0 \simeq 
T_{0,\omega} \times_{\Spec \bK_0} \Spec \bL_0. $$ 
Since  $S_{0,\omega} \times_{\Spec \bK_0} \Spec \bL_0 \simeq S_0 \times_{Y_0} \Spec \bL_0$
and $T_{0,\omega} \times_{\Spec \bK_0} \Spec \bL_0 \simeq T_0 \times_{Y_0} \Spec \bL_0$ there is a
dominant morphism of finite degree $X_0 \to Y_0$ for which $S_0\times_{Y_)} X_0
 \simeq T_0\times_{Y_0}X_0$ and we are done.
\eproof


\section{The Dubouloz-Kishimoto theorem}

The aim of this section is to use the Kraft-Russell Generic Equivalence theorem to get
a rather simple proof of the following
result which has a strong overlap with the Dubouloz-Kishimoto theorem.

\bthm\label{dk.t1} Let $\bk$ be an algebraically closed field of
infinite transcendence degree over the prime field
and let $f : X\to Y$ be a dominant morphism between
geometrically integral algebraic $\bk$-varieties. Suppose that for general closed points 
$y \in Y$, the  fiber $X_y$ contains a Zariski dense subvariety $U_y$ of the form
$U_y \simeq Z_y \times \A_\bk^m$ over a 
$\kappa (y)$-variety $Z_y$. Then there exists a dominant morphism $T\to Y$ 
of finite degree such that
$X_T=X\times_YT$ contains a variety $ W\simeq Z\times \A_\bk^m$ over a $T$-variety $Z$.

\ethm

We start with the following.

\blem\label{dk.l1} Let the notation of  Lemma \ref{ac.l3} hold, $\varphi$ be an 
open immersion,
$\mu$ be a maximal ideal
of the ring $C$ contained in the image of the morphism 
$S_{0,\omega} \times_{\Spec K_0} \Spec  C \to \Spec C$
and $\bk'=C/\mu$ (i.e. the field $\bk'$ is a finite extension of $\bk_0$).
Suppose that $X'=\tX \times_{\Spec C}\Spec \bk'$, $Z'=\tZ \times_{\Spec C}\Spec \bk'$, and
$\varphi' : Z' \to X'$ is the morphism obtained from $\tilde \varphi$ by the base extension
$C \to \bk'$. Let  $Z_0=V_0 \times \A_{\bk_0}^m$ (i.e. $Z=V \times \A_\bk^m$) where $V_0$ is affine.
Then $Z' =V' \times \A_{\bk'}^m$ and $\dim_{\bk'} \varphi' (Z')= \dim_\bk Z$.
\elem

\bproof By the assumption the algebra $\bk_0 [Z_0]$ of regular functions on $Z_0=V_0 \times \A_{\bk_0}^m$
 coincides with   a polynomial ring $B_0[x_1, \ldots , x_m]$ where $B_0$ is some $\bk_0$-algebra.
Hence for $\tZ=Z_0\times_{\Spec \bk_0} \Spec C$ the algebra of regular functions
is $\tB [x_1, \ldots , x_m]$ where $\tB=B_0\otimes_{\bk_0} C$. This implies that 
the algebra $\bk'[Z']$ of regular functions on $Z'$ coincides with $B'[x_1, \ldots , x_m]$ where
$B' =\tB \otimes_C \bk'$ which yields the equality $Z' =V' \times \A_{\bk'}^m$.

Furthermore, by the Noether's normalization lemma there are algebraically independent
elements $y_1, \ldots , y_n \in B_0$ such that $B_0$ is a finitely generated module over
$\bk_0 [y_1, \ldots , y_n]$  and, hence, the natural embedding 
$\iota : \bk_0[x_1, \ldots , x_m,y_1, \ldots , y_n] \to \bk_0 [Z_0]$
is an integral homomorphism (in particular, $\dim Z_0=\dim Z=n+m=:d$). 
Note that $\bk' [Z']= (\bk_0 [Z_0]\otimes_{\bk_0} C)\otimes_C \bk'=\bk_0 [Z_0]\otimes_{\bk_0} \bk'$
and $\iota$ induces a homomorphism $ \bk' [x_1, \ldots , x_n,y_1, \ldots , y_m] \to \bk'  [Z']$
which is integral by \cite[Exercise 5.3]{AM}. Thus $\dim Z'=d$. 
Enlarging $\bk_0$ in this construction
we can suppose that the morphism $\varphi^{-1} : \varphi (Z) \to Z$ can be also obtained
from a morphism of $\bk_0$-varieties. This implies that $\varphi'$ is an immersion and we are done.
\eproof

\noindent {\em Proof of Theorem \ref{dk.t1}}.  Let $\bk_0$, $X_0$, $Y_0$, $\bK_0$ and $X_{0,\omega}$
be as in Notation \ref{ac.n1} and Lemma \ref{ac.l1}. That is,  $\bk_0\subset \bk$
is finitely generated over the prime field, $X = X_0 \times_{\Spec \bk_0} \Spec \bk$,
$Y = Y_0 \times_{\Spec \bk_0} \Spec \bk$, $\bK_0$ is the field of rational functions on $Y_0$,
$X_{0,\omega} = X_0 \times_{\Spec \bk_0} \Spec \bK_0$, and, choosing an
embedding $\bK_0 \hookrightarrow \bk$, 
we have an isomorphism $X_{0,\omega} \times_{\Spec \bK_0}\Spec \bk \to X_y$
for some closed point $y \in Y$.  
Enlarging $\bk_0$ (and, therefore, $\bK_0$)
and treating $\bK_0$ as a subfield of $\bk$ we can  suppose that the
natural immersion $\varphi : Z_y \to X_y$ is obtained from an immersion
$\varphi_0 : Z_y^0 \to X_y^0$ of $\bK_0$-varieties $Z_y^0$ and $X_y^0$
 via the base extension $\Spec \bk \to \Spec \bK_0$.
By Lemma \ref{ac.l3} there exist  a finitely generated $\bK_0$-algebra
$C\subset \bk$, ringed spaces $\tX_y$ and $\tZ_y$ with structure
sheaves consisting of $C$-rings
and a $C$-morphism $\tilde \varphi : \tZ_y \to \tX_y$
such that $X_y=\tX_y\times_{\Spec C} \Spec \bk$, $Z_y=\tZ_y\times_{\Spec C} \Spec \bk$
and $\varphi = \tilde \varphi \times_{\Spec C} \id_{\Spec \bk}$. Let $\mu$ be a maximal ideal
of $C$ and $\bL_0=C/\mu$. By Lemma \ref{dk.l1}
we get an immersion $\varphi' : Z' \to X'$ of the $\bL_0$-varieties 
$X_y'=\tX_y\times_{\Spec C} \Spec \bL_0$ and  $Z_y'=\tZ_y \times_{\Spec C}\Spec \bL_0$
such that $Z' =V' \times \A_{\bL_0}^m$ and $\dim_{\bL_0} \varphi' (Z')=\dim_\bk Z$.

Put $\bL=\bK\otimes_{\bK_0}\bL_0$ where $\bK$ is a field of rational functions on $Y$.
Let $\breve Z = Z' \times_{\Spec \bL_0} \Spec \bL$ (in particular
$\breve Z =\breve V \times \A_\bL^m$), $\breve X = X' \times_{\Spec \bL_0} \Spec \bL$
and $\breve \varphi: {\breve Z} \to \breve X$ be the open immersion induced by $\varphi'$.
By construction the field $\bL_0$ is a finite
extension of $\bK_0$ and, hence, $\bL$ is a finite
extension of $\bK$.  This implies that  there is a dominant morphism $T \to Y$ of finite degree
such that the field of rational functions on $T$ is $\bL$ and
 $X\times_Y T$ is a $\bk$-variety for which the generic fiber of the projection 
 $X\times_Y T$ is $\breve X$. Hence $X\times_Y T$ contains a Zariski open dense subset
 of the form $W \times \A_\bk^m$  which is the desired conclusion. \hfill  $\square$

\brem\label{dk.r1}
 Note that unlike in the original formulation of the Dubouloz-Kishimoto theorem
the field $\bk$ may be countable
and  it may have a positive characteristic.
If $\bk$ has characteristic zero  then similar to Remark \ref{abc.r2}  we can suppose
that $T\to Y$  is  \'etale as in the Dubouloz-Kishimoto theorem.

\erem

\end{document}